\documentclass[11pt, a4paper]{extarticle}

\usepackage[top=90pt, left=90pt, right=90pt, bottom=90pt]{geometry}
\RequirePackage{amsfonts, amsmath, amssymb, amsthm, array, hhline}
\RequirePackage[T2A]{fontenc}
\RequirePackage[utf8]{inputenc}
\RequirePackage{soul}

\usepackage{graphicx, wrapfig}
\graphicspath{{media/}}

\RequirePackage{pgfpages}

\usepackage{hyperref}

\usepackage{esvect}

\def\Id{\mathop{\rm Id}\nolimits}

\def\ort{
\begin{smallmatrix}
\bot \\
$--$
\end{smallmatrix}
}

\newcommand*{\vek}[1]{\vv{#1}}

\newcommand*{\lin}[1]{
\langle
#1 \rangle}

\newcommand*{\otstup}{\vspace{2mm}}

\newtheorem{predl}{\\}[section]

\newcommand{\dok}{{\it Proof.}\hspace{2mm}}
\newcommand{\edok}{$\square$\hspace{2mm}}

\makeindex

\begin{document}

\title{Linear Families of Triangles and Orthology}

\author{Egor Bakaev, Pavel Kozhevnikov$^1$}


\date{
{\footnotesize $^1$ Moscow Institute of Physics and Technology}
\\
\bigskip 
December 21, 2023
}

\maketitle



\newtheorem*{sled}{Сorollary}
\newtheorem*{sled1}{Сorollary 1}
\newtheorem*{sled2}{Сorollary 2}
\newtheorem*{sled3}{Сorollary 3}
\newtheorem*{sled4}{Сorollary 4}
\newtheorem*{sled5}{Сorollary 5}
\newtheorem*{lem}{Lemma}
\newtheorem*{lem1}{Lemma 1}
\newtheorem*{lem2}{Lemma 2}

\newtheorem*{zamech}{Remark}





\begin{figure}
	\centering
	\includegraphics[scale = 1]{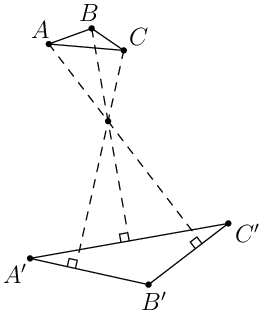}
	\caption{Triangles $ABC$ and $A'B'C'$ are orthologic.}
	\label{fig:orthologic-tr}
\end{figure}

Two triangles are called \textit{orthologic} if the perpendiculars from the vertices of one of them to the sides of the other are concurrent (see Fig.~\ref{fig:orthologic-tr}). 
In this paper, we explore the concept of orthology from various points of view.
Mostly we work in terms of elementary geometry in $\mathbb{R}^2$.
In the final section, we relate the discussed concepts to $\mathbb{R}^4$.
A key idea in this paper is that two orthologic triangles generate a one-parameter (linear) family in which any two triangles are orthologic.
Working with such a family can provide a natural approach to certain questions.

\section{Linear Families of Triangles and Classes of Homothets}

We start with considerations related to linear motion.

\begin{predl}\label{11}
{\bf A single point}
\end{predl}

Let point $A$ move on the plane \textit{linearly}, meaning in a straight line and at a constant speed.
Formally, we define this motion by the mapping $A: \mathbb{R}\to \mathbb{R}^2$, such that the point $A_t = A (t)$ is determined by the radius vector 
$\vek{OA_t} = \vek{OA_0}+ t\cdot \vek{v_a}$, where $O$ is the origin and 
$\vek{v_a}$ is the \textit{velocity vector}.
When $\vek{v_a} \neq \vek{0}$, the \textit{trajectory} $\{A_t\}$ of point $A$ is a specific line $a$ with the direction vector $\vek{v_a}$, and the \textit{graph} of the linear motion of the point is the line $\{(t, A_t)\}$ in $\mathbb{R}^3 = \mathbb{R}\times \mathbb{R}^2$.

\begin{predl}\label{12}
{\bf Pair of points}
\end{predl}

Consider a pair of points $A$ and $B$, each of which moves linearly with velocities $\vek{v_a}$ and $\vek{v_b}$, respectively. Note that $\vek{A_{t}B_{t}} = \vek{A_{0}B_{0}} + t(\vek{v_b}-\vek{v_a})$. In the case when $\vek{A_{0}B_{0}}\parallel (\vek{v_b}-\vek{v_a})$, all vectors $\vek{A_{t}B_{t}}$ are collinear with each other, or equivalently, the graphs of the motions of points $A$ and $B$ in $\mathbb{R}^3$ lie on the same plane. We call this situation and the corresponding pair of moving points $A$ and $B$ \textit{singular}. In this singular case, when $\vek{v_b} \neq \vek{v_a}$, the points collide at some moment (i.e., $A_{t_0} = B_{t_0}$ for some $t = t_0$).

In the \textit{nonsingular} case where $\vek{A_{0}B_{0}}\nparallel (\vek{v_b}-\vek{v_a})$, all vectors in the family $\{\vek{A_{t}B_{t}}\}$ are nonzero. By extending $\vek{A_{\infty}B_{\infty}} \parallel \vek{v_b}-\vek{v_a}$, we can observe that for each direction, there exists a unique vector of the form $\vek{A_{t}B_{t}}$ parallel to that direction (by \textit{direction} we mean a class of parallel lines).

It is easy to see that for a nonsingular pair $A$, $B$ when $ \vek{v_b} \parallel \vek{v_a}$ all lines $A_tB_t$ pass through a single point. 
However, when $ \vek{v_b} \nparallel \vek{v_a}$, it is known that (see, for example, \cite{akopzaslbook}): the envelope for the family of lines $A_tB_t$ is a parabola (at the same time, the lines $a$ and $b$, along which points $A$ and $B$ move, are also lines of the type $A_tB_t$ and are tangent to this parabola). 
Moreover, any circle that passes through the points $A_t$, $B_t$, and $a\cap b$ also passes through the focus of this parabola, denoted as $M_{AB}$.
Thus, $M_{AB}$ is the \textit{Miquel point} for any quadruple from the family of lines $A_tB_t$. Additionally, $M_{AB}$ is the center of spiral similarities that match pairs of vectors of the form $\vek{A_tB_t}$. Furthermore, $M_{AB}$ is the center of the spiral similarity that transforms line $a$ into $b$ and matches points $A_t$ and $B_t$ for a fixed $t$.

In $\mathbb{R}^3$, the graphs of linear motion of a nonsingular pair consist of
a pair of skew lines, and lines of the form $(t, A_tB_t)$ constitute a family of straight-line generatrices of a hyperbolic paraboloid.

\begin{predl}\label{13}
{\bf Triple of points} 
\end{predl}

Next, we consider a triple of points $A$, $B$, and $C$, each of which moves linearly with velocities $\vek{v_a}$, $\vek{v_b}$, and $\vek{v_c}$ respectively. 
We call such a triple \textit{nonsingular} if every pair of points within it is nonsingular. 
The triple of linearly moving points $A, B, C$ defines a set of triangles $\{A_tB_tC_t\}$ ($t\in \mathbb{R}$), which we refer to as a \textit{linear family of triangles}. 
Usually, we consider \textit{nonsingular linear families of triangles}, i.e., the nonsingular case of the linear motion of three points.

In the nonsingular case, in any triangle $A_tB_tC_t$, none of the three vertices coincide. However, $A_t, B_t, C_t$ may lie on a straight line; in this case, we call $A_tB_tC_t$ a \textit{degenerate triangle}.


It is clear that a linear family can be generated by any two of its triangles (say, taking them as $A_0B_0C_0$ and $A_1B_1C_1$), and any two triangles of a linear family determine the entire family.


Let us mention some other known geometric facts related to the nonsingular linear family $\{A_tB_tC_t\}$.
The centers of spiral similarities $M_{AB}$, $M_{BC}$, and $M_{AC}$ (defined in \ref{12}) lie on the so-called similarity circle of three shapes (as shapes, one can simply take the vectors $\vek{A_0A_1}$, $\vek{B_0B_1}$, and $\vek{C_0C_1}$). 
The line connecting $a \cap b$ with $M_{AB}$, and two analogous lines intersect at a single point located on the similarity circle. 
The common point $P_t$ of the three circles — namely, the circle passing through $A_t$, $B_t$, $a \cap b$, and two analogous circles — moves along the similarity circle. 
These facts can be found, for example, in \cite[section <<The similarity circle of three figures>>]{pras}. (These facts have become quite popular and have been included in many problem books and sheets on spiral similarity.)

\bigskip

\begin{predl}\label{14}
{\bf Triple of points $\leftrightarrow$ pair of vectors   } 
\end{predl}

The triangle $ABC$ can be defined, up to a translation (shift), by an ordered pair of (Euclidean) vectors $(\vek{b}, \vek{c}) = (\vek{AB}, \vek{AC})$.
For the linear family $A_tB_tC_t$, we define $\vek{b}_t = \vek{A_{t}B_{t}}$ and $\vek{c}_t = \vek{A_{t}C_{t}}$.
It is straightforward to observe the following vector equalities:
\begin{equation}\label{linlin}
\vek{b}_t= (1-t) \vek{b}_0 +t \vek{b}_1, \,\,\,\,\,
\vek{c}_t= (1-t) \vek{c}_0 +t \vek{c}_1,
\end{equation}
from which, in particular, it is clear that  replacing triangles $A_0B_0C_0$ and $A_1B_1C_1$ with certain shifts of them leaves the family $\{A_tB_tC_t\}$ unchanged, meaning it consists of the same triangles up to a shift.

\begin{predl}\label{16}
{\bf Classes of homothets}
\end{predl}

Next, we examine equivalence in relation to homotheties and translations. 
We say that two sets $S_1$ and $S_2$ on the plane are equivalent if $S_2$ is derived from $S_1$ by a certain translation or homothety. 
The equivalence class $[S]$ of the set $S$ consists of homothets and translations of the set $S$; this class can be succinctly termed as the {\it class of homothets} of the set $S$.
For a nondegenerate triangle $ABC$, its class $[ABC]$ is uniquely determined by an ordered triple of directions parallel to $BC$, $CA$, and $AB$.

Equivalence of triangles $ABC$ and $A'B'C'$ corresponds to the proportionality of pairs of vectors $(\vek{AB}, \vek{AC})$ and $(\vek{A'B'}, \vek{A'C'})$. Thus, the class $[ABC]$ corresponds to the set of pairs of the form $(t\vek{AB}, t\vek{AC})$.
Given a nonsingular linear family of triangles $A_tB_tC_t$, we define a {\it linear family of classes} $\{ [A_tB_tC_t] \}$, supplemented by the class of the <<triangle at infinity>> $[A_{\infty}B_{\infty}C_{\infty}]$, which can be defined by a pair of vectors $(\vek{A_{\infty}B_{\infty}}, \vek{A_{\infty}C_{\infty}}) $, proportional to the pair $(\vek{v_b}-\vek{v_a}, \vek{v_c}-\vek{v_a})$.

The right-hand sides of the equality (\ref{linlin}), written in <<homogeneous form>>
\begin{equation}
\label{linob}
( x_0\vek{b}_0+x_1\vek{b}_1, \,\,\, x_0 \vek{c}_0+x_1\vek{c}_1 ),
\end{equation}
where $x_0, x_1$ are arbitrary constants, which are not simultaneously equal to 0, describe a linear family of classes generated by the classes $[A_{0}B_{0}C_{0}]$ and $[A_{1}B_{1}C_{1}]$. The case $x_0+x_1=0$ corresponds to the class $[A_{\infty}B_{\infty}C_{\infty}]$ of the <<triangle at infinity>>.

From the formulas (\ref{linob}), it is easy to see that the linear family of classes remains unchanged by replacing triangles $A_0B_0C_0$ and $A_1B_1C_1$ by their homothets. Moreover, a linear family of classes is determined by any two of its (distinct) classes. Therefore, the following is true: two non-coinciding linear families of classes either do not intersect or intersect at exactly one class.

\begin{predl}\label{15}
{\bf Degenerate triangles}
\end{predl}

A well-known statement, which was posed as the <<pedestrians problem>> in the All-Soviet Union Olympiad (see \cite[problem 228]{vaseg}), is as follows:

{\it In a linear family of triangles, 
either there are no more than two degenerate triangles, or all the triangles of the family are degenerate. }


\textit{Proof.} This statement can be proven, for example, using the equality (\ref{linlin}).
Suppose in the linear family, there are two degenerate triangles. Without loss of generality, we  assume these are triangles $A_0B_0C_0$ and $A_1B_1C_1$.
Hence $\vek{b}_0 \parallel \vek{c}_0$ and $\vek{b}_1 \parallel \vek{c}_1$.
If the factors of proportionality  are the same, that is, $\vek{b}_t = \lambda \vek{c}_t$ for $t=0$ and $t=1$, then according to (\ref{linlin}) 
the equality holds for all $t$. 
(In this case, we can say that the entire family consists of degenerate triangles that are <<similar>> to each other.)
Otherwise, it is clear that 
$\vek{b}_t \nparallel \vek{c}_t$ when $t\neq 0$ and $t\neq 1$.
\edok

If all triangles of the family are degenerate, then the family is called {\it degenerate}. For a degenerate nonsingular family of triangles, 
the motion graphs of its vertices belong to one hyperbolic paraboloid.

\begin{predl}\label{17}
{\bf Examples}
\end{predl}

{\bf (1) <<Spiral-similar>> family}

\begin{figure}
	\centering
	\includegraphics[scale = 1]{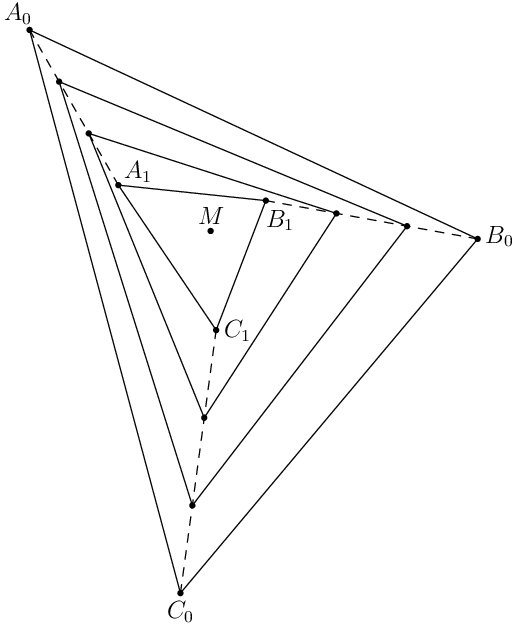}
	\caption{Four similar triangles of a single <<spiral-similar>> family.}
	\label{fig:similar-tr}
\end{figure}

Consider a nonsingular linear family generated by triangles $A_0B_0C_0$ and $A_1B_1C_1$, which are similar and have the same orientation (in the nondegenerate case). From the equalities (\ref{linlin}), it follows that any triangle $A_tB_tC_t$ in the family is also similar to $A_0B_0C_0$ and has the same orientation.

Moreover, for all triangles in the family, there exists a common corresponding point $M$. This point serves as the center of the spiral similarity that aligns any pair of vectors of the form $\vek{A_tB_t}$, and therefore any pair of triangles in our family. Thus, in the considered specific case, using notations of \ref{12} and \ref{13}, we have $M_{AB} = M_{BC} = M_{CA} = M$.

\otstup

{\bf (2) <<Pedal>> family }

Consider fixed lines $a, b, c$, with a point $P$ that moves linearly. The projections $A_t, B_t, C_t$ of the point $P_t$ onto the lines $a, b, c$ move linearly as well, forming a linear family of {\it pedal} triangles. If the lines $a, b, c$ form a triangle $\Delta$, then the pedal triangle of the point $P$ may degenerate (into the Simpson line), and this occurs if and only if $P$ lies on the circumcircle of $\Delta$.

The construction of the <<pedal>> family can be generalized by considering, instead of the pedal triangle of point $P_t$, points $A_t$, $B_t$, $C_t$ on the lines $a, b, c$ such that each of the oriented angles $\angle (P_tA_t, a)$, $\angle (P_tB_t, b)$, $\angle (P_tC_t, c)$ equals a fixed value $\alpha$, so that the directions of the lines $P_tA_t$, $P_tB_t$, $P_tC_t$ are fixed; in this case, the triangle $A_tB_tC_t$ can be termed as an <<$\alpha$-pedal>>.

In the case of the $\alpha$-pedal family, 
the similarity circle (along which $P_t$ moves) degenerates into a line, 
and the points $M_{AB}$, $M_{BC}$, $M_{CA}$ become fixed points on this line. 
(Consequently, the lines connecting the pairs of points $M_{AB}$ and $a\cap b$, 
$M_{BC}$ and $b\cap c$, $M_{CA}$ and $c\cap a$, are parallel).

Conversely, for collinear points $M_{AB}$, $M_{BC}$, $M_{CA}$, there corresponds an 
<<$\alpha$-pedal>> family: one can move the point $P_t$ linearly along the line 
$M_{AB}M_{BC}M_{CA}$. Then the point $A_t$, which is the intersection of circles passing through the triples of points $P_t$, $a\cap c$, $M_{CA}$ and $P_t$, $a\cap b$, $M_{AB}$, moves linearly along the line $a$, 
and the line $A_tP_t$ retains its direction.

\begin{predl}\label{155}
{\bf Geometry of $\{T_t\}$ in the case of two degenerate triangles}
\end{predl}

Let there be two degenerate triangles in the linear nonsingular family of triangles $A_tB_tC_t$, lying on lines $x$ and $y$, intersecting at point $O$. For simplicity, we denote the triangle $A_tB_tC_t$ as $T_t$. 
We reparameterize the family so that the vertices of the triangle $T_0$ are on the line $y$, 
and the vertices of the triangle $T_1$ are on the line $x$. 
The conditions $f(A_0)= A_1$, $f(B_0)= B_1$, and $f(C_0)= C_1$ define a unique mapping $f: y\to x$, 
which preserves the cross ratios of points (here $x$ and $y$ are considered as extended by points at infinity).
For a fixed constant $\lambda$, on each line of the form $P f(P)$ (where $P$ runs through line~$y$), 
we mark the point $P_{\lambda}$ such that $\vek{PP_{\lambda}} = \lambda \vek{Pf(P)}$, 
in particular, $P_0 = P$ and $P_1=f(P)$. (Note that these notations align with the notation $A_{\lambda}$, 
$B_{\lambda}$, and $C_{\lambda}$ for the vertices of the family triangle.)

\otstup

{\it Case 1.} In the case where the ratios are equal, $\vek{A_0B_0}/\vek{A_0C_0}= \vek{A_1B_1}/\vek{A_1C_1}$ 
(in other words, when $A_0B_0C_0$ and $A_1B_1C_1$ are degenerate similar triangles),
in accordance with \ref{15}, all triangles in the family $\{T_t\}$ 
are degenerate and similar to each other. In the considered case, $f$ is linear (i.e., preserves ratios),
and according to \ref{12}, lines of the form $Pf(P)$ are tangent to a fixed parabola 
(in particular, lines $x$ and $y$ are tangent to it).
For a fixed $\lambda$, the points $P_{\lambda}$ lie on a line, also tangent to the same parabola.

\otstup

{\it Case 2.} Now suppose $\vek{A_0B_0}/\vek{A_0C_0}\neq \vek{A_1B_1}/\vek{A_1C_1}$. The following holds:
{\it For a fixed $\lambda \neq 0, 1$, points of the form $P_{\lambda}$ form a hyperbola $\gamma_{\lambda}$ with asymptotes parallel to lines $x$ and~$y$.}

{\it Proof.} 
In the affine coordinate system $Oxy$ (with coordinate axes $x$ and $y$), 
one can specify the coordinates of point $P$ as $(0, \tau)$ and consequently, the coordinates of 
$f(P)$ as $(\alpha (\tau), 0)$, where $\alpha$ is a linear fractional function (since $f$ 
preserves double ratios). 
(Note that in this case, $\alpha$ is not a linear function, as this corresponds to the case 
$\vek{A_0B_0}/\vek{A_0C_0}= \vek{A_1B_1}/\vek{A_1C_1}$.)
Thus, $P_{\lambda}$ has coordinates $(\lambda \alpha(\tau) , (1-\lambda) \tau)$. 
Therefore, $P_{\lambda}$ forms a hyperbola $x= \lambda \alpha(\frac{y}{1-\lambda})$.
\edok

\otstup
Note that the hyperbola $\gamma _{\lambda}$ is uniquely defined by the following conditions: 
it contains the vertices $A_{\lambda}, B_{\lambda}, C_{\lambda}$ and its asymptotes are parallel to lines $x$ and $y$.

For $\lambda \neq 0, 1$ and $\mu \neq 0, 1$, there exists 
a unique affine transformation $g_{\lambda, \mu}$ of the plane that maps the vertices 
$A_{\lambda}, B_{\lambda}, C_{\lambda}$ to $A_{\mu}, B_{\mu}, C_{\mu}$ respectively. 
We see that in the $Oxy$ coordinates, the transformation $g_{\lambda, \mu}$ can be defined as $$x^* = \frac{\mu}{\lambda} x, \ \  
y^* = \frac{1-\mu}{1-\lambda} y,$$ where $(x^*, y^*)$ is the image of point $(x, y)$. 
In this transformation, each point $P_{\lambda}$ maps to $P_{\mu}$, and consequently,
the hyperbola $\gamma _{\lambda}$ maps to the hyperbola $\gamma _{\mu}$.

\begin{figure}
	\centering
	\includegraphics[scale = 1]{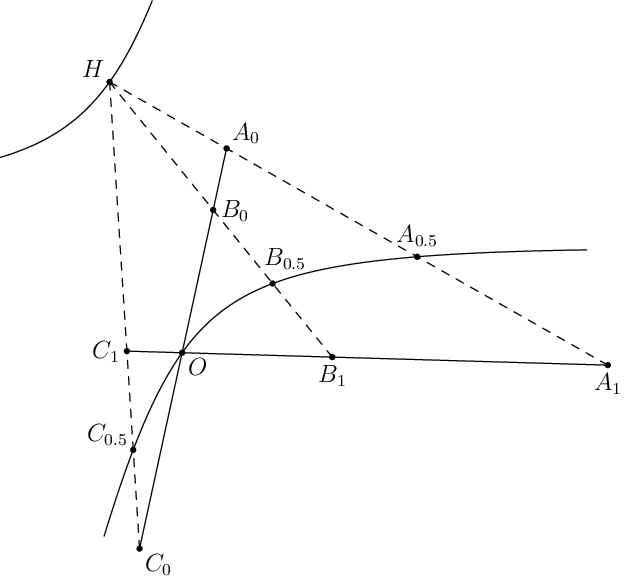}
	\caption{Case 2.1. Hyperbola $\gamma_{\lambda}$ for $\lambda = 0.5$.}
	\label{fig:2}
\end{figure}

\otstup
{\it Case 2.1.} Consider the {\it concurrent} case where the lines $a$, $b$, $c$ (along which the vertices 
of the triangles $T_t$ move) intersect at a single point $H$ (see Fig. \ref{fig:2}). 
Then the mapping $f: y\to x$ coincides with the central projection with center $H$, hence all lines 
of the form $Pf(P)$ pass through $H$.
Note that for a fixed $\lambda \neq 0, 1$ 
the point $H$ is of the form $P_{\lambda}$ for some point $P$,
and also $O = P_{\lambda}$ for $P=O$, therefore
{\it the hyperbolas $\gamma_{\lambda}$ pass through $H$ and $O$.}

\otstup 

This fact can also be established through direct substitution into the equation: we find $\alpha (\tau) = 
\frac{x_0}{1-\frac{y_0}{\tau}}$, where $(x_0, y_0)$ are the coordinates of the point $H$, from which
the equation for $\gamma_{\lambda}$ is $x= \frac{\lambda x_0 y}{y-(1-\lambda)y_0}$.
\edok

Furthermore, $\gamma_{\lambda}$ has asymptotes 
$y=(1-\lambda) y_0$ and $x=\lambda x_0$ 
(the geometric position of the asymptotes can be described 
by considering the directions of the lines $Pf(P)$, close to the direction of the lines $x$ and $y$).
Note that all hyperbolas $\gamma_{\lambda}$ have a pair of common points and parallel asymptotes, hence they belong to a single pencil.

\otstup

{\it Case 2.2.} In the {\it nonconcurrent} case, it is known that lines of the form $Pf(P)$ 
(when $P$ runs along the line $y$, and $f(P)$ along the line $x$, preserving cross-ratios) are tangent to a certain conic 
$\varepsilon$, in particular, the lines $x$ and $y$ are 
tangents to $\varepsilon$. In fact, $\varepsilon$ is either an ellipse or a hyperbola, 
and the parabola case arises above in 1).
Note that the hyperbolas $\gamma_{\lambda}$ are tangent to the conic $\varepsilon$ twice 
(the tangency of $\gamma_{\lambda}$ and $\varepsilon$ corresponds 
to the tangency of the line $Pf(P)$ and $\varepsilon$ at the point $P_{\lambda}$).

\begin{predl}\label{1555}
{\bf Case of two degenerate triangles, one of which is $T_{\infty}$}
\end{predl}

Let there be exactly two degenerate triangles in the linear nonsingular family $\{T_t\}$,
the first of which lies on the line $y$,
and the second is $T_{\infty}$.

According to \ref{16}, the case of the degenerate triangle $T_{\infty}$ corresponds to
the parallelism of the velocity differences $(\vek{v_b}-\vek{v_a})\parallel (\vek{v_c}-\vek{v_a})$.
In other words, in the affine coordinate system $Oxy$ with the $Ox$ axis parallel to the vector $(\vek{v_b}-\vek{v_a})$,
the velocity vectors $\vek{v_a}$, $\vek{v_b}$, $\vek{v_c}$ have equal ordinates.

Next, we introduce an affine coordinate system, taking the line $y$ as the ordinate axis,
and drawing the abscissa axis $Ox$ parallel to the vector $(\vek{v_b}-\vek{v_a})$.
We reparameterize the family so that the vertices of the triangle $T_0$ lie on the line $y$.
As a result, $A_0 = (0, \tau_a)$, $B_0 = (0, \tau_b)$, $C_0 = (0, \tau_c)$.
Moreover, assume that the (equal) ordinates of the vectors $\vek{A_0A_1}$, $\vek{B_0B_1}$, $\vek{C_0C_1}$
are equal to 1, so that $\vek{A_0A_1} = (\xi_a , 1 )$, $\vek{B_0B_1} = (\xi_b , 1 )$, 
$\vek{C_0C_1} = (\xi_c , 1 )$ (where the abscissas $\xi_a$, $\xi_b$, $\xi_c$ are different, otherwise
$\{T_t\}$ is singular).
Then $A_{\lambda}$ has coordinates $(\lambda \xi_a, \tau_a + \lambda)$.
The coordinates of the points $B_{\lambda}$ and $C_{\lambda}$ can be written similarly.

The conditions $\alpha(\tau_a)= \xi_a $, $\alpha(\tau_b)= \xi_b $, and $\alpha(\tau_c)= \xi_c $ determine a unique 
linear fractional function~$\alpha$. 
Correspondingly, for each point $P=P_0 (0, \tau) $ on the line $y$, we define  
$P_{\lambda}  = (\lambda \alpha (\tau), \tau + \lambda) $
(which is consistent with the coordinates of $A_{\lambda}$,  $B_{\lambda}$, $C_{\lambda}$).
Note that the mapping that associates point $P$ with the line $\{P_{\lambda}\}$, 
actually defines a unique
mapping $f$ from the line $y$ to the line at infinity, 
preserving cross ratios such that $f(A_0) = a$, 
$f(B_0) = b$, $f(C_0) = c$.
The following statement holds: 

\otstup

{\it 
For a fixed $\lambda \neq 0$, points of the form $P_{\lambda}$ form a hyperbola $\gamma_{\lambda}$ 
with asymptotes parallel to $x$ and $y$.
}

{\it Proof.} 
From the coordinates $P_{\lambda}  = (\lambda \alpha (\tau), \tau + \lambda) $ 
in our affine coordinate system $Oxy$, we see that
$P_{\lambda}$ moves along the hyperbola defined by $x= \lambda \alpha(y-\lambda)$.
(It is easy to notice that the case when $\alpha$ turns out to be a linear function corresponds to the collinearity
of the points $A_1$, $B_1$, $C_1$, which contradicts the conditions.)
\edok

\otstup

Note that the hyperbola $\gamma_{\lambda}$ is uniquely determined by the following conditions: 
it contains the vertices $A_{\lambda}, B_{\lambda}, C_{\lambda}$ and its asymptotes are parallel to the lines $x$ and $y$.

For $\lambda \neq 0$ and $\mu \neq 0$, there exists 
a unique affine transformation $g_{\lambda, \mu}$ of the plane that maps the vertices 
$A_{\lambda}, B_{\lambda}, C_{\lambda}$  to $A_{\mu}, B_{\mu}, C_{\mu}$ respectively. 
We see that in the $Oxy$ coordinates, the transformation $g_{\lambda, \mu}$ can be defined as 
$x^* = \frac{\mu}{\lambda} x$, 
$y^* = y + \mu - \lambda$. 
In this case, each point $P_{\lambda}$ is mapped to  $P_{\mu}$, and correspondingly
the hyperbola $\gamma_{\lambda}$ is transformed into the hyperbola $\gamma_{\mu}$.

\otstup

{\it Case 0.1.} Consider the {\it concurrent case}, where the lines $a$, $b$, $c$ intersect at a single point $H(x_0, y_0)$. 
Then the aforementioned mapping $f$ 
coincides with the central projection with the center at $H$, hence all the lines
of the form $\{P_{\lambda}\}$ (with fixed $P$) also pass through $H$.
We show that 

\otstup

{\it the hyperbolas $\gamma_{\lambda}$ pass through $H$ and have the asymptote $y$.}

{\it Proof.} Indeed, let us write down the condition that for a given $P$, the line $\{P_{\lambda}\}$ passes through the point $H(x_0, y_0)$:
$\lambda_0 \alpha (\tau) = x_0$, $\tau + \lambda_0 = y_0$ holds for some $\lambda_0$.
By eliminating $\lambda_0$, we get $\alpha (\tau) = \frac{x_0}{y_0-\tau}$ (for all $\tau$). 

Therefore, for a fixed $\lambda$, the hyperbola $\gamma_{\lambda}$  
has the equation $x = \frac{\lambda x_0}{\lambda + y_0 - y}$. As can be seen, $x=0$ is an asymptote, 
and $(x_0, y_0)$ satisfies the equation.
\edok

\otstup

{\it Case 0.2.} In the {\it nonconcurrent} case, lines of the form $\{P_{\lambda}\} $  (for a fixed $P$) are tangent to some 
parabola $\varepsilon$. In particular, it is tangent to the line $y$ and the line at infinity. 
Note that, as in \ref{155}, the hyperbolas $\gamma_{\lambda}$ touch the parabola $\varepsilon$ twice.

\section{Orthology Relation}

\begin{predl}\label{21}
{\bf Definition of orthology} 
\end{predl}

Recall that two triangles $ABC$ and $A'B'C'$, where $A'B'C'$ is nondegenerate 
(and $ABC$ possibly degenerate),
are {\it orthologic}, if the perpendiculars from $A$ to $B'C'$, from $B$ to $C'A'$,
and from $C$ to $A'B'$, intersect at a single point. 
The fact of orthology of triangles $ABC$ and $A'B'C'$ we denote $ABC\ort A'B'C'$. 
If we label the triangles $ABC$ and $A'B'C'$ (with a fixed order of vertices) as $T$ and $T'$, 
then we get a shorter notation for orthology: $T\ort T'$.
The common intersection point of the perpendiculars is the
{\it center} of orthology; denoted $O_{T, T'}$.

\begin{predl}\label{22}
{\bf Reflexivity}
\end{predl}

Clearly, orthology is reflexive: $T\ort T$,
and the center of orthology $O_{T, T}$ is the orthocenter of triangle~$T$.

\begin{predl}\label{23}
{\bf Symmetry}
\end{predl}

\begin{figure}
	\centering
	\includegraphics[scale = 1]{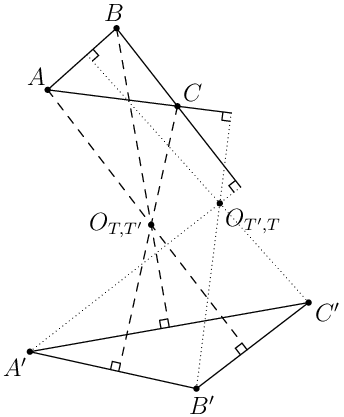}
	\caption{Triangles $ABC$ and $A'B'C'$ are orthologic. Points $O_{T', T}$ and $O_{T, T'}$ are their centers of orthology.}
	\label{fig:orthologic-tr2}
\end{figure}

It is known that the condition for orthology for nondegenerate triangles is symmetric, i.e., 
$T\ort T'$ $\Rightarrow$ $T'\ort T$. Hence, if $T\ort T'$, then
there exists a second center of orthology 
$O_{T', T}$, generally distinct from $O_{T, T'}$ (see Fig.\ref{fig:orthologic-tr2}).

The symmetry of the orthology relation follows, for instance, from the well-known Carnot's theorem
\begin{equation}\label{carnot}
A'B^2-A'C^2+B'C^2-B'A^2+C'A^2-C'B^2=0,
\end{equation}
which is equivalent to the condition for orthology $ABC\ort A'B'C'$. 

\otstup
It is clear that if one moves the vertex of one triangle along the perpendicular to the corresponding side 
of another triangle (say, move point $A$ along the perpendicular to $B'C'$), then
both the condition for orthology and the left side of
(\ref{carnot}) remain unchanged.

\otstup
We take the equality (\ref{carnot}) as the definition of orthology for a pair of degenerate triangles.

\begin{predl}\label{24}
{\bf Orthologic quadruples}
\end{predl}

Let $D = O_{T, T'}$ and $D' = O_{T', T}$ be the centers of orthology for nondegenerate 
orthologic triangles $ABC$ and $A'B'C'$. 
We have $AB\perp C'D'$, $AC\perp B'D'$, etc. — six perpendiculars for quadruples of 
points $A, B, C, D$ and $A', B', C', D'$. As we see, in this construction, the four pairs of points $A$ and $A'$, 
$B$ and $B'$, $C$ and $C'$, $D$ and $D'$ are equivalent. Therefore, for instance, 
$A$ and $A'$ are centers of orthology for orthologic triangles $BCD$ and $B'C'D'$, and so on.
We can call the quadruples $A, B, C, D$ and $A', B', C', D'$ orthologic.
The six aforementioned perpendiculars are taken as the definition of orthologic tetrahedra 
$ABCD$ and $A'B'C'D'$ in space (see, for instance, \cite{turgor}),
in this sense, orthologic quadruples on a plane can be considered degenerate orthologic tetrahedra.
It is clear (for instance, from the symmetry of the orthology relation) that 
from the five aforementioned perpendicularities, the sixth follows.

We prove that {\it orthologic quadruples are affinely equivalent}, in other words, 
an affine transformation that maps $A, B, C$ respectively into $A', B', C'$, 
maps $D = O_{T, T'}$ to $D' = O_{T', T}$. This fact is known as the Rideau's theorem 
(see, for example, \cite{zasl_sondat}).

\dok
Let $\vek{DA} = \vek{a}$, $\vek{DB} = \vek{b}$, $\vek{DC} = \vek{c}$,
$\vek{D'A'} = \vek{a'}$, $\vek{D'B'} = \vek{b'}$, $\vek{D'C'} = \vek{c'}$.
The case when among the vectors $\vek{a}$, $\vek{b}$, $\vek{c}$,  
$\vek{a'}$, $\vek{b'}$, $\vek{c'}$ there exists a zero vector is easily dealt with. 
Henceforth, we assume all these vectors are nonzero.

From the perpendicularities, we deduce $\vek{a'}(\vek{b}-\vek{c}) = 0$, hence
the dot products $\vek{a'} \vek{b}$ and $\vek{a'} \vek{c}$ are equal. 
All six analogous dot products equal the same value, which we denote $p$.
Also note that $\vek{a} \vek{a'} \neq p$, otherwise from the equalities
$\vek{a} \vek{a'} = \vek{b} \vek{a'} = \vek{c} \vek{a'}=p$, the collinearity 
of points $A$, $B$, $C$ would follow.

Let $\alpha   \vek{a} + \beta \vek{b} + \gamma   \vek{c} = \vek{0}$ and 
$\alpha '  \vek{a'} + \beta' \vek{b'} + \gamma '  \vek{c'} = \vek{0}$,
where $\alpha  + \beta  + \gamma   = \alpha ' + \beta'  + \gamma ' = 1$.
To complete the proof and establish the required affine equivalence, 
it suffices to show that $\alpha = \alpha'$, $\beta = \beta'$, $\gamma = \gamma'$.
By scalar multiplying the first equality by $\vek{a'}$, we get 
$\alpha \vek{a} \cdot \vek{a'} + (\beta  + \gamma)p = 0$, from which 
$\alpha (\vek{a} \cdot \vek{a'}-p)  + (\alpha + \beta  + \gamma)p = 0$, 
from which $\alpha$ is uniquely determined as $\alpha = \frac{p}{p-\vek{a} \vek{a'}}$.
By similar reasoning, we find that $\alpha '$ equals the same value.
\edok

We also present a generalization of Rideau's theorem, discovered by A. Myakishev (\cite[problem 4.2]{turgor-zamechat}): if the pairs of triangles $A_iA_jA_k$ and $B_iB_jB_k$ are orthologic for all $\{i,j,k\}\subset\{1,2,3,4\}$, then the quadrilaterals $A_1A_2A_3A_4$ and $B_1B_2B_3B_4$ are affinely equivalent.

\bigskip 
\bigskip 

\begin{predl}\label{25}
{\bf Orthology as a relation on classes of homothets}
\end{predl}

From the definition, it is clear that the fact $T \ort T'$ remains valid when replacing triangle $T'$ with its translation or a homothet. Due to the symmetry of orthology, the same holds true for $T$. Thus, the orthology relation between triangles $T \ort T'$ elevates to the relation $[T] \ort [T']$ on classes of homothets or, equivalently, on (ordered) triples of directions.

Clearly, if a translation or homothety is applied to one of the two orthologic triangles, say $ABC$, the entire quadruple $A, B, C, O_{T, T'}$ undergoes this translation or homothety (meanwhile, the quadruple $A', B', C', O_{T', T}$ remains unchanged).

\begin{predl}\label{29}
{\bf Orthology as a relation on pairs of vectors}
\end{predl}

Since the fact of orthology $ABC \ort A'B'C'$ is retained when one of the triangles is translated, orthology can be understood as a relation on a set of ordered pairs of vectors (as before, for the triangle $ABC$, we associate the pair of vectors $(\vek{b}, \vek{c}) = (\vek{AB}, \vek{AC})$).

The relation (\ref{carnot}) can be transformed into many equivalent forms using the dot product. In particular, (\ref{carnot}) is equivalent to the equality $\vek{AB}\cdot \vek{A'C'} - \vek{AC}\cdot \vek{A'B'}=0$.

Thus, the orthology relation on pairs of vectors $(\vek{b}, \vek{c}) \ort (\vek{b'}, \vek{c'})$ can be expressed by the equation
\begin{equation}\label{ort_pary}
 \vek{b} \cdot \vek{c'} - \vek{b'}  \cdot \vek{c}=0.
\end{equation}

\begin{predl}\label{210}
{\bf Orthology of degenerate triangles}
\end{predl}

Let $ABC$ and $A'B'C'$ be degenerate triangles, such that $\vek{c} = \alpha \vek{b}$ and $\vek{c'} = \alpha ' \vek{b'}$. According to (\ref{ort_pary}), 
$$ABC \ort A'B'C'  
\ \ \Longleftrightarrow \ \ 
\vek{b}\cdot (\alpha ' \vek{b'}) - (\alpha  \vek{b}) \cdot \vek{b'} = 0 
\ \ \Longleftrightarrow \ \ 
(\alpha -\alpha ') \vek{b} \cdot \vek{b'} = 0.$$ 
We see that $ABC \ort A'B'C'$ holds in two cases: when $\vek{b} \perp \vek{b'}$, i.e., when lines $ABC$ and $A'B'C'$ are perpendicular, and also in the case $\alpha =\alpha '$, in other words, when the degenerate triangles $ABC$ and $A'B'C'$ are <<similar>>.

\otstup

\begin{predl}\label{2111}
{\bf The case of coinciding orthology centers} 
\end{predl}

\begin{figure}
	\centering
	\includegraphics[scale = 1]{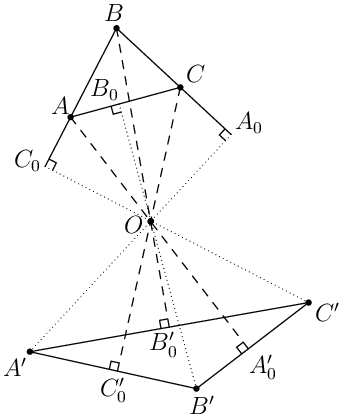}
	\caption{The case of coinciding orthology centers.}
	\label{fig:orthologic-tr3}
\end{figure}

Let $ABC$ and $A'B'C'$ be two orthologic triangles with coinciding orthology centers (see Fig.\ref{fig:orthologic-tr3}): 
$O=O_{T, T'}= O_{T', T}$, then $OA\perp B'C'$, $OB\perp C'A'$, $OC\perp A'B'$,  
$OA'\perp BC$, $OB'\perp CA$, $OC'\perp AB$. 
Let $C_0 = AB\cap OC'$, $C'_0 = A'B'\cap OC$, and so on, such that $A_0B_0C_0$ is the pedal triangle
of the point $O$ relative to the triangle $ABC$ and analogously,
$A'_0B'_0C'_0$ is the pedal triangle of the point $O$ relative to $A'B'C'$. 
Then $OC_0B$ and $OB'_0C'$ are similar right-angled triangles, 
from which we derive $\vek{OB'_0}\cdot \vek{OB} = \vek{OC_0}\cdot \vek{OC'}$. All analogous 
products are also equal. 
This situation can be described, for example, by stating that $A'B'C'$ and $ABC$ are respectively the images of the pedal triangles $A'_0B'_0C'_0$ and $A_0B_0C_0$ under an inversion (possibly of an imaginary radius) centered at $O$.
Or equivalently, 
the lines containing the sides of one of the triangles $T$, $T'$, 
are the polars of the corresponding vertices of the other triangle with respect to the circle centered at $O$. 

It is known (see, for example, problem 20 in \cite{akopzaslbook}), that 
{\it in the case $O_{T, T'}= O_{T', T}$ the triangles $T$ and $T'$ are perspective
(i.e., $AA'$, $BB'$, and $CC'$ are concurrent).}

\otstup

According to \ref{25}, from the general case of orthologic triangles $ABC$ and $A'B'C'$, one can 
obtain the case of coinciding orthology centers 
by shifting one of the quadruples $A, B, C, O_{T, T'}$ and
$A', B', C', O_{T', T}$.

\begin{predl}\label{27}
{\bf Orthology and Maxwell's theorem} 
\end{predl}

Another relation can be introduced (it can be called {\it harmonicity}) on the classes of homothets of nondegenerate triangles (or on triples of distinct directions). 
We say that the class $[ABC]$ is {\it harmonic} to the class $[A'B'C']$ if the lines passing through $A, B, C$ and 
parallel to the lines $B'C'$, $C'A'$, $A'B'$ respectively, intersect at a single point.

Let the class $[ABC]$ correspond to the triple of directions $a, b, c$, 
and the class $[A'B'C']$ to the triple of directions $a', b', c'$. Subsequently,
the orthology of the triple $(a, b, c)$ to the triple $(a', b', c')$ is equivalent to 
the harmonicity of the triple $(a, b, c)$ to the triple $(a'^{\bot}, b'^{\bot}, c'^{\bot})$ 
(where $x^{\bot}$ denotes the direction perpendicular to $x$). 
The latter is equivalent to the harmonicity of the triple $(a^{\bot}, b^{\bot}, c^{\bot})$ to the triple $(a', b', c')$, 
since the simultaneous rotation of all directions by the same angle preserves both harmonicity and orthology. 
We see that the facts about the symmetry of the relations of orthology and harmonicity are equivalent. 
The fact about symmetry of the harmonicity relation 
is known as Maxwell's theorem \cite{max}.

\otstup 

Next, we present \textit{a scheme for an {\it affine} proof} of Maxwell's theorem 
without using orthology.
Suppose that the lines $a, b, c, a', b', c'$, parallel respectively to the sides $BC$, $CA$, $AB$, 
$B'C'$, $C'A'$, $A'B'$ of the triangles $ABC$ and $A'B'C'$, pass through a single point $O$. Let us intersect these lines with an arbitrary line $l$ and obtain intersection points 
$A_0, B_0, C_0, A'_0, B'_0, C'_0$ respectively. The condition of harmonicity of $ABC$ and $A'B'C'$ is equivalent to the equality
\begin{equation}\label{harm}
\dfrac{\vek{A_0B'_0}}{\vek{B'_0C_0}}\cdot \dfrac{\vek{C_0A'_0}}{\vek{A'_0B_0}}
\cdot \dfrac{\vek{B_0C'_0}}{\vek{C'_0A_0}} = -1
\end{equation} 
(this can be proven, for example, by using Ceva's theorem in its sine form),
from which symmetry of relation is apparent. \edok

\otstup
The term {\it harmonicity} was chosen by us because the equality (\ref{harm})
has the form <<the cyclic ratio of $2n$ points is equal to $-1$>>,
and, as is known, a similar equality for $n=2$ defines a harmonic quadruple.

Consider two quadruples $A, B, C, D$ and $A', B', C', D'$, in which there are no three collinear points.
Let $AB\parallel C'D'$, $AC\parallel B'D'$, $AD\parallel B'C'$,  
$BC\parallel A'D'$, $BD\parallel A'C'$, $CD\parallel A'B'$ (so that $ABC$ and $A'B'C'$ are harmonic,
or, equivalently, $ABD$ and $A'B'D'$ are harmonic, etc.).
Then the quadruples $A, B, C, D$ and $A', B', C', D'$ are affinely equivalent.
This fact is equivalent to Rideau's theorem on the affine equivalence of orthologic quadruples (but of course
it can also be proven within the framework of affine geometry, without using orthology).

\begin{predl}\label{28}
{\bf $\alpha$-orthology}
\end{predl}

We say that the triangle $ABC$ is $\alpha$-{\it orthologic} to the triangle $A'B'C'$ if 
there exists a point $P$ such that $$\angle (AP, B'C') = \angle (BP, C'A') = 
\angle (CP, A'B') = \alpha.$$ In this case, <<ordinary>> orthology is
$\frac{\pi}{2}$-orthology, and 
the previously introduced relation of harmonicity coincides with $0$-orthology.

If the triangle $ABC$ is $\alpha$-orthologic to the triangle $A'B'C'$,
then after rotation by $\alpha$ (or by $\alpha + \frac{\pi}{2}$)
it becomes respectively harmonic (orthologic) to the triangle $A'B'C'$.

Hence, it is clear that if $ABC$ is $\alpha$-orthologic to the triangle $A'B'C'$, 
then $A'B'C'$ is $(-\alpha)$-orthologic to the triangle $ABC$.

\begin{predl}\label{211}
{\bf Examples}
\end{predl}

\otstup

{\bf (1) Pedal triangle}

\begin{figure}
	\centering
	\includegraphics[scale = 1]{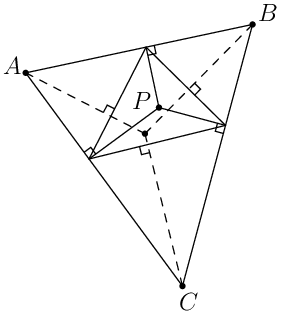}
	\caption{Pedal triangle of point $P$.}
	\label{fig:pedal}
\end{figure}

If $T$ is the pedal triangle of point $P$ for the triangle $ABC$ (see Fig.\ref{fig:pedal}),
then $T \ort ABC$, and the orthology centers are
the point $P$ and its isogonal conjugate point.

More generally, the $\alpha$-pedal triangle of point $P$ is $\alpha$-orthologic to the triangle $ABC$
(with the same centers of $\alpha$-orthology and $(-\alpha)$-orthology).

\otstup

{\bf (2) Radical axes}

The radical axis of two circles with centers $O_1$ and $O_2$ and radii $R_1$ and $R_2$ is the set of points $X$ such that $XO_1^2 - XO_2^2 = R_1^2 - R_2^2$. The radical axis is a line perpendicular to
$O_1O_2$. In the case of intersecting circles, the radical axis passes through the points of intersection of the circles.

Let $\omega_i$, for $i=1, 2, 3$, be three circles with centers $O_i$, not lying on a single line.
If $X_{ij}$ is an arbitrary point on the radical axis of circles $\omega_i$ and $\omega_j$,
then $O_1O_2O_3\ort X_{23}X_{31}X_{12}$ (for example, it is easily verified by (\ref{carnot})). This fact is equivalent to the statement that
the three radical axes intersect at one point (the radical center).

\otstup

{\bf (3) Replacement by homothet}

\begin{figure}
	\centering
	\includegraphics[scale = 1]{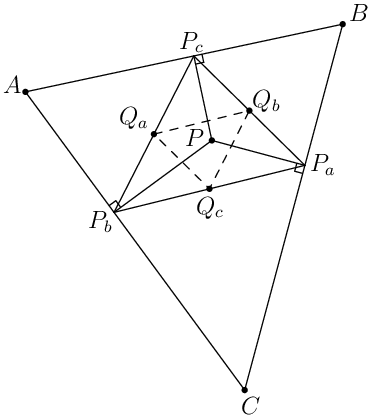}
	\caption{For the problem of replacement by homothet.}
	\label{fig:pedal2}
\end{figure}

The idea that orthology remains unchanged
when one of the triangles is replaced by its homothet is very useful.
Here we illustrate it with the solution of the following problem.

{\it Problem.} Consider the projections $P_a, P_b, P_c$ of the point $P$ onto the lines containing the sides of triangle $ABC$ (Fig. \ref{fig:pedal2}).
Through the midpoint of $P_bP_c$, draw a line $l_a$ perpendicular to $BC$.
Similarly, lines $l_b$ and $l_c$ are defined. Prove that $l_a, l_b, l_c$ are concurrent.

{\it Solution.}
Let $Q_a$, $Q_b$, $Q_c$ be the midpoints of the corresponding sides of the triangle
$P_aP_bP_c$. The statement of the problem means that $Q_aQ_bQ_c \ort ABC$. But this is true since
$P_aP_bP_c\ort ABC$, and $P_aP_bP_c$ and $Q_aQ_bQ_c$ are homothetic. \qed

\section{Orthology and Linear Families}

Let us continue to explore linear families of triangles.

\begin{predl}\label{31}
{\bf Orthology to a given triangle. Linearity of $O_{T_t,  T'}$}
\end{predl}

{\it
Suppose two triangles from the linear family $\{T_t\}$ (i.e., triangles $T_t$ at two different values of $t$) are orthologic to a given triangle $T'$. Then $T_t \ort T'$ for any $t$.

Moreover, if $T'$ is nondegenerate, then the center of orthology $O_{T_t,  T'}$ moves linearly.
}

\dok
First, suppose $T'$ is nondegenerate.
Assume that $T_{\lambda}\ort T'$ and $T_{\mu}\ort T'$ for $\lambda\neq \mu$.
Let $O_t$ be a point moving linearly such that $O_{\lambda} = O_{T_{\lambda},  T'}$ and 
$O_{\mu} = O_{T_{\mu},  T'}$.
Since $A_{\lambda}O_{\lambda} \perp B'C'$ and $A_{\mu}O_{\mu} \perp B'C'$, we obtain that
$A_tO_t \perp B'C'$ for any $t$. By making similar arguments for the other vertices, we conclude that
the perpendiculars from $A_t$ to $B'C'$, from $B_t$ to $C'A'$, and from $C_t$ to $A'B'$ intersect
at point $O_t$, i.e., $O_t = O_{T_t,  T'}$ for any $t$.

To prove the first statement without using the non-degeneracy of $T'$,
it is sufficient to use the linearity of the dot product and, taking into account (\ref{linlin}) and (\ref{ort_pary}),
derive the equality $\vek{b'}\vek{c}_t - \vek{c'}\vek{b}_t = 0$ for any $t$ from the equalities $\vek{b'}\vek{c}_0 - \vek{c'}\vek{b}_0 = 0$ and $\vek{b'}\vek{c}_1 - \vek{c'}\vek{b}_1 = 0$.
\edok

\otstup 

Thus, for a nondegenerate triangle $T'$, the trajectory of $O_{T_t, T'}$ is a line or degenerates to a point.
It is easy to see that degeneracy occurs in the case when the velocity vectors of the vertices of the triangle $T_t$
are perpendicular to the corresponding sides of the triangle $T'$. In particular, if the triangle $T'$
belongs to the family $\{T_t\}$, then the trajectory of $O_{T_t, T'}$ can degenerate to a point --- the orthocenter
of the triangle~$T'$.

Note also that the trajectory of $O_{T_t, T'}$ does not change when replacing $T'$ with its homothet or shift.

\begin{predl}\label{32}
{\bf Orthologic linear family}
\end{predl}

We briefly refer to a linear family of triangles in which any two triangles are orthologic as an {\it orthologic family} of triangles. Similarly, we define orthologic families of classes of homothets.

The following is a sufficient condition for a family to be orthologic.

{\it If some two triangles from the family $\{T_t\}$ are orthologic to each other, then $\{T_t\}$ is an orthologic family.}

\dok
Without loss of generality, assume that $T_0\ort T_1$.
Since $T_0\ort T_0$, it follows from \ref{31} that
$T_0\ort T_t$ for all $t$. Similarly, $T_1\ort T_t$ for all $t$.
Now, for each fixed $\lambda$, it was proved that $T_1\ort T_{\lambda}$ and $T_0\ort T_{\lambda}$.
From this, it again follows from \ref{31} that $T_t\ort T_{\lambda}$ for all $t$.
\edok

\otstup

Similarly, {\it if in a linear family of classes of homothets some two different classes, say $[T_0]$ and $[T_1]$,
are orthologic, then this family is orthologic.}

Moreover, in such a case (i.e., when $[T_0]\neq [T_1]$), it turns out that
{\it a triangle $T$ is orthologic to both triangles $T_0$ and $T_1$ if and only if
$[T]$ belongs to the linear family generated by $[T_0]$ and $[T_1]$,}
thus orthologic families of classes are maximal sets in which any pair of triangles are orthologic.
The explanation of this fact is provided in \ref{57}, though a proof can also be obtained through direct calculation.

\begin{predl}\label{34}
{\bf Degenerate triangles of orthologic families}
\end{predl}

If there are two degenerate triangles in an orthologic family, then, according to \ref{210},
either these triangles are similar, or they lie on perpendicular lines.
In the first case, according to \ref{15}, the family is degenerate, consisting of degenerate <<similar>> triangles.

Thus, in a nondegenerate orthologic family, two degenerate triangles can only lie on perpendicular lines.
It turns out that this statement can be strengthened to the following.

{\it If $\{[T_t]\}$ is a nonsingular nondegenerate orthologic linear  family of classes, then there are exactly two classes of degenerate triangles in it, and these classes correspond to perpendicular directions.}

The proof of this fact can be found below in \ref{55},
though it can also be proven directly.

\otstup

Thus, a nonsingular nondegenerate orthologic family $\{T_t\}$ contains exactly two degenerate triangles,
one of which may be the <<triangle at infinity>> $T_{\infty}$.

In the first case, similar to \ref{155},
we reparameterize the orthologic family
so that the degenerate triangles are $T_1$ and $T_0$,
and they lie respectively on the axes $Ox$ and $Oy$
of the 
Cartesian coordinate system $Oxy$.

In the case of degeneration of $T_{\infty}$,
we assume that the other degenerate triangle is $T_0$, and it lies on the line $y$.
This case corresponds to the conditions $\vek{v_a}-\vek{v_b}\perp y$ and $\vek{v_a}-\vek{v_c}\perp y$.
In other words, $T_{\infty}$ degenerates if the projections of the velocities of vertices $A_t$, $B_t$, $C_t$ on the axis $y$ are equal,
i.e., the projection of $A_tB_tC_t$ onto $y$ is a shift of the degenerate triangle $A_0B_0C_0$ (along $y$).

\begin{predl}\label{36}
{\bf Examples of the use of orthologic families}
\end{predl}

{\bf (1) Orthopole}

Let $T_0$ denote the orthogonal projection of the triangle $T_1$ onto a specific line $\ell$.
As is known (and easily verified by~(\ref{carnot})),
 $T_0 \ort T_1$. The center of orthology $O_{T_0, T_1}$
is called the {\it orthopole} of the triangle $T_1$ and the line~$\ell$.
The orthologic linear family generated by $T_0$ and  $T_1$ contains 
the images of $T_1$ under all possible dilations, contractions, and reflections with respect to $\ell$.
In particular, from this example, we see that any triangle is orthologic to its mirror reflection
(and therefore, to any shift of the mirror reflection).
Of course, the orthology of a triangle and its reflection can be proved directly using (\ref{ort_pary}) or by explicitly specifying the position of one of the centers of orthology (it lies on the circumcircle of the triangle).

In accordance with \ref{34}, in the considered family, $T_{\infty}$ is degenerate.

\otstup

{\bf (2) Criterion for the orthology of pedal triangles}

Let $ABC$ be a given triangle, $O$ be the center of its circumcircle,
$A_0B_0C_0$ and $A_1B_1C_1$ be the pedal triangles of points $P_0$ and $P_1$ respectively.

We establish the following criterion, which was also presented as Problem 17 in the correspondence round of the Sharygin Olympiad in 2009:
$A_0B_0C_0\ort A_1B_1C_1$ $\Leftrightarrow$
$P_0, P_1, O$ are collinear.

The triangles $A_0B_0C_0$ and $A_1B_1C_1$ generate a family of pedal triangles $A_tB_tC_t$ of points
$P_t$, moving along the line $P_0P_1$. Since each of the triangles in this family is orthologic to $ABC$,
the family $A_tB_tC_t$ is orthologic if and only if the corresponding family
of classes $[A_tB_tC_t]$ contains the class $[ABC]$.
But the pedal triangle is homothetic to $ABC$ only for the point $O$.

Another explanation is related to the consideration of degenerate triangles of the family $A_tB_tC_t$ ---
the Simson lines of the intersection points of the line $P_0P_1$ with the circle $(ABC)$.
The family $A_tB_tC_t$ is orthologic if and only if there are two such lines that are perpendicular; that is,
when the points of intersection of the line $P_0P_1$ with the circle $(ABC)$ exist and
are diametrically opposed.

\otstup

{\bf (3) L. Emelyanov's Problem}

{{\it Problem.} (All-Russian Olympiad, 2002.) Prove that the lines passing through the points of tangency of the exscribed circles with the sides of the triangle, parallel to the corresponding bisectors, intersect at one point.}

{\it Solution.} Let $X_a, X_b, X_c$ be the points of tangency 
of the sides of the triangle $ABC$ with the exscribed circles. 
It is enough to show that the triangle $X_a X_b X_c$ is orthologic to some triangle $T$, whose
sides have the directions of the external bisectors of the triangle $ABC$.

{\it First way.} As $T$, we take the triangle $K_aK_bK_c$ with vertices at the points 
of tangency with the inscribed circle.
Since $K_a$ and $X_a$ are symmetric with respect to the midpoint $A_0$ of the side $BC$ (and a similar statement 
is true for the pairs $K_b$, $X_b$ and $K_c$, $X_c$), the linear family generated by triangles 
$K_aK_bK_c$ and $X_aX_bX_c$ contains the medial triangle $A_0B_0C_0$. Since 
$K_aK_bK_c \ort ABC$, then $K_aK_bK_c \ort A_0B_0C_0$ (replacement by homothet). 
Therefore, our family is orthologic, from which $K_aK_bK_c\ort X_aX_bX_c$.

{\it Second way.} In fact, this problem can be reduced to Example 2, as $K_aK_bK_c$ is the pedal triangle 
of the center $I$ of the inscribed circle, and $X_aX_bX_c$ is the pedal triangle of the Bevan point, which is symmetric to $I$ with respect to the circumcenter.

{\it Third way.} (Proposed by A. Buchaev.)
Let $T$ be the triangle $A'B'C'$, where $A'$, $B'$, and $C'$ are the midpoints of the arcs $BAC$, $CBA$, and $ACB$, respectively.
It is easy to show that $BX_c=CX_b$, which implies that the triangles $A'X_cB$ and $A'X_bC$ are congruent.
This equality results in $A'X_c=A'X_b$. Given this and similar equalities,
the orthology $A'B'C' \ort X_aX_bX_c$ follows immediately from (\ref{carnot}).
\edok

\otstup

{\bf (4) The problem about the midpoints of the altitudes}

{\it Problem.} (European Mathematical Cup, 2013.)
Let $X$, $Y$, and $Z$ be the midpoints of the altitudes $AD$, $BE$, and $CF$ of the triangle $ABC$
respectively. Prove that the perpendiculars from $D$ to $YZ$, from
$E$ to $ZX$, and from $F$ to $XY$, intersect at one point.

{\it Solution.} 
Since $X$, $Y$, and $Z$ are located on the altitudes of triangle $ABC$, we have 
$ABC\ort XYZ$. Consequently, $ABC$ and $XYZ$ generate an orthologic linear family.
This family also contains the triangle $DEF$ (if $ABC$ corresponds to the value $t=0$,
and $XYZ$ corresponds to the value $t=1$, then $DEF$ corresponds to the value $t=2$).
Since in this family any two triangles
are orthologic, we have $DEF\ort XYZ$, which is the required statement.
\edok

\otstup
\bigskip

{\bf (5) Proof of Sondat's theorem}

\begin{figure}
	\centering
	\includegraphics[scale = 1]{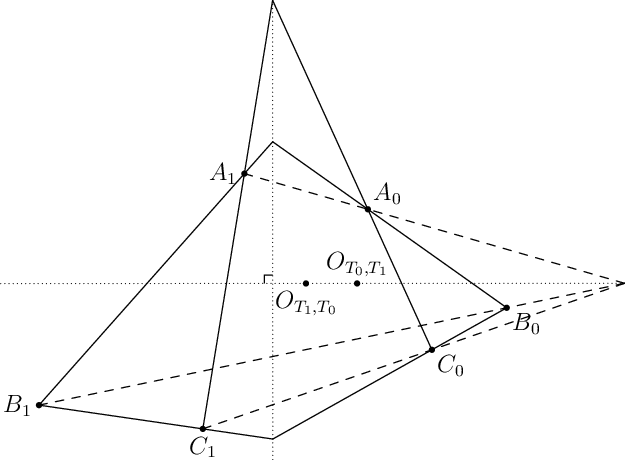}
	\caption{Illustration for Sondat's theorem.}
	\label{fig:sondat}
\end{figure}

Let $A_0B_0C_0$ and $A_1B_1C_1$ be two given orthologic and perspective triangles.
Sondat's theorem asserts that then the centers of orthology and perspector lie on a straight line,
and this line is perpendicular to the axis of perspective (Desargues' line) of the given triangles (Fig. \ref{fig:sondat}).

\dok
The triangles $A_0B_0C_0$ and $A_1B_1C_1$ generate an orthologic family $A_tB_tC_t$.
Let $H$ be the perspector of the triangles $A_0B_0C_0$ and $A_1B_1C_1$, which is also the perspector of any two triangles in the family. 
It suffices to show that for $t\neq 0$, the vector $\vek{HO_{T_t, T_0}}$ is perpendicular to the axis of perspective of the triangles $A_0B_0C_0$ and $A_tB_tC_t$.
Hence, for $t=1$, we have that 
$\vek{HO_{T_1, T_0}}$ is perpendicular to the axis of perspective of the triangles 
$A_0B_0C_0$ and $A_1B_1C_1$. By swapping the roles of $T_0$ and $T_1$, 
we also prove that $\vek{HO_{T_0, T_1}}$ is perpendicular to the axis of perspective of the triangles 
$A_0B_0C_0$ and $A_1B_1C_1$, and from this the statement of Sondat's theorem will follow.

Note that as the point $O_{T_t, T_0}$ moves linearly, 
the vector $\vek{u}(t) = \vek{HO_{T_t, T_0}}$ depends linearly on $t$ (that is, the coordinates of this vector 
in any Cartesian coordinate system are linear functions of $t$). 

Next, we demonstrate that the point $A_0B_0\cap A_tB_t$ moves linearly (along the line $A_0B_0$).
It is similarly established that the point $A_0C_0\cap A_tC_t$ moves linearly.
Subsequently, the direction vector of the axis of perspective, $\vek{v}(t)$, connecting these points,
also depends linearly on $t$.

To prove this, we select a coordinate system in which the line $A_0B_0$ coincides with the $y$-axis,
so that the point $A_0B_0\cap A_tB_t$ has coordinates $(0, y_1)$,
and points $A_t$ and $B_t$ have coordinates $(\alpha t, \beta t + \gamma)$ and $(\alpha ' t, \beta' t + \gamma')$
respectively, for certain constants $\alpha$, $\beta$, $\gamma$, $\alpha'$, $\beta'$, and $\gamma'$.
The collinearity condition of the vectors connecting $A_t$ and $B_t$ with $A_0B_0\cap A_tB_t$ is written as
$$\alpha t (\beta' t + \gamma' - y_1) = \alpha' t (\beta t + \gamma - y_1).$$ 
Upon eliminating $t$ from this relation, it becomes clear that $y_1$ is linearly dependent on~$t$.

Thus, our objective is to prove that $\vek{u}(t) \perp \vek{v}(t)$ for all $t$, or equivalently, $(\vek{u}(t), \vek{v}(t))\equiv 0$, where the vectors $\vek{u}(t)$ and $\vek{v}(t)$ depend linearly on $t$.
Since the dot product $(\vek{u}(t), \vek{v}(t))$ is a polynomial in $t$ of degree no more than 2, it suffices 
to identify three different values $t_i$, for $i=1, 2, 3$, satisfying $\vek{u}(t_i) \perp \vek{v}(t_i)$. 
Let $t_1$ correspond to the position $A_{t_1}=H$. 
Consequently, the axis of perspective
of triangles $A_0B_0C_0$ and $A_{t_1}B_{t_1}C_{t_1}$ coincides with the line $B_0C_0$,
and the vector $\vek{HO_{T_{t_1}, T_0}} = \vek{A_{t_1}O_{T_{t_1}, T_0}}$ is perpendicular to $B_0C_0$, as required.
By a similar argument, let $t_2$ and $t_3$ correspond to the positions $B_{t_2}=H$ and $C_{t_3}=H$.
In the case of a nonsingular family $A_tB_tC_t$, the three parameter values $t_1$, $t_2$, $t_3$ are pairwise different, since the points $A_t$, $B_t$, $C_t$ are pairwise different, and in this case the theorem is proven.

If the family $A_tB_tC_t$ is singular, then the triangles of this family has one of its sides retaining its orientation, 
$A_0B_0\parallel A_tB_t$ for any $t$. Then in the reasoning above, we can make the following simplifications: the vector $\vek{v}(t)$ retains its direction parallel to $A_0B_0$, and we need to establish that
the vector $\vek{u}(t)$ is always perpendicular to $A_0B_0$. Since $\vek{u}(t)$ depends linearly on $t$, it is enough
to verify this for two different parameter values. Suitable are $t=t_2$ and $t=t_3$,
corresponding to the positions $B_{t_2}=H$ and $C_{t_3}=H$, except in the case where $A_{t_2}=B_{t_2}=C_{t_2}=H$
and the family $A_tB_tC_t$ consists of homothetic triangles. But in the latter case, the theorem is obvious.
\edok

\otstup

{\bf (6) The problem of flies on altitudes}

Here we present the problem by E. Bakaev, which actually served as the starting point for writing this paper.

\otstup

{\it Problem.} Three flies landed on different vertices of triangle $ABC$ and crawled along the lines containing its altitudes at constant speeds. At some point, all flies found themselves on the same line $x$. After some more time, they all ended up on another line $y$. Prove that $x \perp y$.

\otstup

{\it Solution 1.}
The movement of the flies defines an orthologic linear family of triangles $A_tB_tC_t$. Indeed, $A_0B_0C_0 = ABC$, and clearly, $A_tB_tC_t \ort ABC$ for any $t$. Therefore, the statement of the problem follows immediately from \ref{34}.
\edok

\otstup

We also offer an author's elementary geometric solution, which additionally proves in the conditions of the problem that in the case of a nonsingular family $A_tB_tC_t$ there exist lines $x$ and $y$ (i.e., degenerate triangles of the family $A_tB_tC_t$).

\otstup

{\it Solution 2.}
Assume that $ABC$ is $A_0B_0C_0$.

Let the perpendiculars to the altitudes $a, b, c$ of triangle $A_0B_0C_0$, drawn respectively through points $A_t$, $B_t$, $C_t$, form a triangle $\Delta_t$. It follows that $A_tB_tC_t$ is the pedal triangle of point $H$ for triangle $\Delta_t$. Note that $\{ \Delta_t \}$ is a linear (singular) family of homothetic (or aligned by translation) triangles. These triangles have a common center $P$ of homothety, which is at the same time a degenerate triangle of the family $\{ \Delta_t \}$ (in the case of parallel translation, $P$ can be considered a point at infinity).

Therefore, all the circles $\Omega_t$, circumscribed around triangles $\Delta_t$, are homothetic with the center at $P$, and the center $S_t$ of the circle $\Omega_t$ moves linearly along the line passing through $P$. Note that $H=S_0$, hence $S_t$ moves along the line $PH$. (Observe that $P \neq H$, otherwise all triangles $A_tB_tC_t$ are homothetic with the center at $H$.)

The degeneracy of the pedal triangle $A_tB_tC_t$ (into the Simson line) is equivalent to the inclusion of $H$ in the circle $\Omega_t$. From the homothety with the center $P$, it follows that $H \in \Omega_t$ for two different values of $t=t_1$ and $t=t_2$. In these positions, the homothety transforming $\Delta_{t_1}$ into $\Delta_{t_2}$ sends the point $H$ on the circle $\Omega_{t_1}$ to the point diametrically opposed to $H$ on the circle $\Omega_{t_2}$. As is known, the Simson lines of a triangle corresponding to diametrically opposed points are perpendicular (and homothety preserves the directions of lines). Therefore, the degenerate triangles $A_{t_1}B_{t_1}C_{t_1}$ and $A_{t_2}B_{t_2}C_{t_2}$ lie on perpendicular lines.
\edok

\begin{predl}\label{55_2}
{\bf The trajectory of $O_{T', T_t}$ is a conic}
\end{predl}

{\it
Let $\{T_t\}$ be a nonsingular linear family of triangles, orthologic to a given triangle $T'$. It follows that the orthology center $O_{T', T_t}$ moves along a conic passing through the vertices $A'$, $B'$, $C'$ of the triangle~$T'$.
}

\otstup

\dok
 The coordinates of the vectors of the sides of the triangle $T_t$ depend linearly on $t$. Rotating these vectors by $90^{\circ}$, we obtain the direction vectors $\vek{u}_a (t)$, $\vek{u}_b(t)$, $\vek{u}_c(t)$ of the perpendiculars $A'O_{T', T_t}$, $B'O_{T', T_t}$, $C'O_{T', T_t}$, drawn respectively through $A'$, $B'$, $C'$ to $B_tC_t$, $C_tA_t$, $A_tB_t$. These vectors also depend linearly on $t$. Therefore, there exists an affine transformation $f$ of the plane, mapping $A'$ to $B'$ and $\vek{u}_a(t)$ to $\vek{u}_b(t)$, such that $f$ maps the line $A'O_{T', T_t}$ to $B'O_{T', T_t}$ (for any $t$).

As is known (see, for example, \cite{akopzaslbook}), if $f$ is a projective transformation mapping a pencil of lines $\Pi(X)$ (with center $X$) to a pencil $\Pi(Y)$ (where $Y\neq X$), then the points of intersection $l\cap f(l)$, where $l\in \Pi(X)$, lie on a fixed line or on a conic passing through $X$ and $Y$. In our case, it follows that $O_{T', T_t}$ moves along a line or along a conic passing through the points $A'$ and $B'$ (and similarly, through $C'$).
\edok

\otstup
It is clear that the set $\{O_{T', T_t}\}$ remains unchanged when replacing the family $\{T_t\}$ with another family that defines the same family of classes $\{[T_t]\}$. When the triangle $T'$ is translated or homothetically transformed, the set $\{O_{T', T_t}\}$ undergoes the corresponding translation or homothetic transformation.

If the triangle $T'$ is degenerate, then $\{O_{T', T_t}\}$ can be considered a degenerate conic, which is the union of the line $A'B'C'$ with the line along which the point $O_{T', T_t}$ moves.

Note that the asymptotic directions of the conic $\{O_{T', T_t} \}$ are perpendicular to the lines containing the degenerate triangles of the family $\{T_t\}$. Therefore, from \ref{34}, it follows:

\otstup

{\it
A linear nonsingular nondegenerate family $\{T_t\}$ is orthologic if and only if the conic of orthology centers ${O_{T', T_t}}$ has a pair of perpendicular asymptotic directions (i.e., it is a rectangular hyperbola or a pair of perpendicular lines).
}

\otstup

Let $\{ T_t\}$ be a nondegenerate nonsingular orthologic family, and let $T_{\lambda}$ be its fixed nondegenerate triangle. 
It follows from the above discussion that $\{O_{T_{\lambda}, T_t}\}$ is a hyperbola with asymptotes parallel to the lines containing the degenerate triangles of the family. Moreover, this hyperbola obviously contains the points $A_{\lambda}, B_{\lambda}, C_{\lambda}$. Thus, $\{O_{T_{\lambda}, T_t}\}$ satisfies the conditions that uniquely define the hyperbola $\gamma_{\lambda}$ from the considerations in \ref{155} and \ref{1555}, that is


\begin{center}
	{\it $\{O_{T_{\lambda}, T_t}\}$ coincides with $\gamma_{\lambda}$.}
\end{center}


{\bf Example: image of a line under isogonal conjugation}

Let there be a triangle $ABC$, and let $\{T_t\}$ be the pedal family of triangles (see \ref{17}) for a point $P_t$ linearly moving along a line $p$. It follows that $ABC$ is orthologic to each of the triangles $T_t$, and in accordance with  Example 1 from \ref{211}, the conic $\{ O_{ABC, T_t} \}$ is the image of the line $p$ under isogonal conjugation with respect to $ABC$.

From Example 2 in \ref{36}, we know that the pedal family $T_t$ is orthologic if and only if $p$ passes through $O$, or equivalently, when $\{ O_{ABC, T_t} \}$ passes through the orthocenter $H$ of the triangle $ABC$. In this case, $\{ O_{ABC, T_t} \}$ is a rectangular hyperbola or a pair of perpendicular lines, consistent with \ref{55_2}.


For example, if $p$ is the line $OL$, where $L$ is the intersection of the symmedians (the Lemoine point), then $\{ O_{ABC, T_t} \}$ is the Kiepert hyperbola. The same hyperbola can be obtained as the set of orthology centers for another family of triangles (which, of course, defines the same family of homothet classes). This family $\{A_tB_tC_t\}$ is generated by the medial triangle $A_0B_0C_0$ of the triangle $ABC$ and the triangle $A_1B_1C_1$ with vertices at the centers of the squares constructed on the sides of the triangle $ABC$ outside of it. Thus, $A_t$, $B_t$, $C_t$ are the vertices of similar isosceles triangles constructed on the sides of the triangle $ABC$ either outside or inside it. It is known that the set of orthology centers $O_{ABC, A_tB_tC_t}$ is the Kiepert hyperbola. The Kiepert hyperbola is also the set of perspectors of $ABC$ and $A_tB_tC_t$. It is easy to check   that $AA_{1/t}\perp B_tC_t$, hence the perspector of $ABC$ and $A_{1/t}B_{1/t}C_{1/t}$ is the orthology center $O_{ABC, A_tB_tC_t}$, and vice versa. Also note that the points $O_{ABC, A_tB_tC_t}$ and $O_{ABC, A_{1/t}B_{1/t}C_{1/t}}$ lie on a line with $O$ (see \cite{zaslkie}).

\begin{predl}\label{371}
{\bf The line connecting the orthology centers $O_{T_{\lambda}, T_{\mu}}$ and $O_{T_{\mu}, T_{\lambda}}$}
\end{predl}

Let $\{T_t\}$ be a nondegenerate nonsingular orthologic family, and let $T_{\lambda}$ and $T_{\mu}$ be two of its nondegenerate triangles.

As noted in \ref{55_2}, the orthology centers $O_{T_{\lambda}, T_{\mu}}$ and $O_{T_{\mu}, T_{\lambda}}$ lie on the corresponding hyperbolas $\gamma_{\lambda}$ and $\gamma_{\mu}$, defined in \ref{155} and \ref{1555}. Moreover, the affine transformation $g_{\lambda, \mu}$ that maps $A_{\lambda}$, $B_{\lambda}$, $C_{\lambda}$ to $A_{\mu}$, $B_{\mu}$, $C_{\mu}$, respectively, according to Rideau's theorem (see \ref{24}) maps $O_{T_{\lambda}, T_{\mu}}$ to $O_{T_{\mu}, T_{\lambda}}$. But it follows, in the notation of \ref{155} and \ref{1555}, that the points $O_{T_{\lambda}, T_{\mu}}$ and $O_{T_{\mu}, T_{\lambda}}$ are respectively $P_{\lambda}$ and $P_{\mu}$ (for the same point $P$ on the line $y$).  From this, the following corollaries arise:

\otstup

1) Let $\{T_t\}$ be a concurrent family. Then the line $P_0P_1$ (which is also $P_{\lambda}P_{\mu}$) passes through the perspector $H$. Thus, we obtain that

\begin{center}
{\it $O_{T_{\lambda}, T_{\mu}}$, $O_{T_{\mu}, T_{\lambda}}$, and $H$ are collinear.}
\end{center}

In fact, we have obtained another proof of the collinearity statement from Sondat's theorem.
\otstup 

2) Let $\{T_t\}$ be a nonconcurrent family. From \ref{155} and \ref{1555}, we know that the lines $P_0P_1$ are tangent to the conic $\varepsilon$ (i.e., the same conic that the lines $a, b, c$, as well as the lines containing degenerate triangles of the family, are tangent to). Therefore,
\begin{center}
{\it the line $O_{T_{\lambda}, T_{\mu}} O_{T_{\mu}, T_{\lambda}}$ is tangent to the conic $\varepsilon$.}
\end{center}

In some sense, this statement can be considered a generalization of the collinearity statement from Sondat's theorem.

Next, note that

\begin{center}
{\it the line $\{O_{T_t, T_{\lambda}}\}$ is also tangent to the conic $\varepsilon$.}
\end{center}
	
Indeed, let the line $\{O_{T_t, T_{\lambda}}\}$ (for a fixed $\lambda$) and the hyperbola $\{O_{T_{\lambda}, T_t}\}$, which have a common point $O_{T_{\lambda}, T_{\lambda}}$, intersect once more at a point $X$. Then $X$ is of the form $O_{T_{\lambda}, T_r}$. 
But the point $O_{T_r, T_{\lambda}}$ also lies on the line $\{O_{T_t, T_{\lambda}}\}$. Thus, the line $\{O_{T_t, T_{\lambda}}\}$ coincides with the line $O_{T_r, T_{\lambda}} O_{T_{\lambda}, T_r}$, and the required statement follows from the previous one.

\begin{predl}\label{371000}
{\bf Additional remarks on a concurrent orthologic family}
\end{predl}

We continue to consider a nondegenerate nonsingular orthologic family of triangles $\{T_t\}$.   In this section, we assume that the points $A_t$, $B_t$, $C_t$ move linearly along three different lines $a$, $b$, $c$, which intersect at the point $H$, so that any two triangles in our family are perspective with the perspector $H$.

\otstup

{\bf Perspector as an orthocenter}

Next, we prove that
\otstup

{\it For the family $\{T_t\}$, there exists a unique $h$ (possibly, $h=\infty$), such that $A_hB_h\perp c$, $B_hC_h\perp a$, $C_hA_h\perp b$.}
\otstup

In other words, if $h\neq \infty$, then there is a triangle in the family for which $H$ is the orthocenter. In this case, the perspective orthologic family is described by example 6 from \ref{36}.

\otstup 

\dok
If $A_{\infty}B_{\infty}\perp c$, $B_{\infty}C_{\infty}\perp a$, $C_{\infty}A_{\infty}\perp b$, then the proof is complete.

Otherwise, suppose for definiteness that $A_{\infty}B_{\infty}$ is not perpendicular to $c$. From \ref{12}, it follows that there is exactly one real $h$ for which $A_hB_h\perp c$. If the orthocenter of triangle $A_hB_hC_h$ coincides with $H$, then the proof is complete. Otherwise,  if it does not and is instead located at a different point $H'\neq H$
 (clearly lying on line $c$), let us consider the point $H'' = O_{T_s, T_h}$ for some $s\neq h$. Since $C_sH''\perp A_hB_h$, it follows that $H''\in c$. Consequently, there is a homothety with center at $H$ which maps triangle $A_sH''B_s$ onto triangle $A_hH'B_h$. From this, $A_sB_s\parallel A_hB_h$, which contradicts the assumption that $\{T_t\}$ is nonsingular.
\edok

Thus, for the found parameter $h$, it is true that $H = O_{T_{\lambda}, T_h}$ for any $\lambda$.
In particular, we obtain another explanation of the fact that
the hyperbolas $\{O_{T_{\lambda}, T_t} \}$ pass through $H$.

\otstup

{\bf Coincidence of the orthology centers $O_{T_{\lambda}, T_{\mu}} = O_{T_{\mu}, T_{\lambda}}$ at the point $O$}


Let our family contain two degenerate triangles $T_0$ and $T_1$ lying on perpendicular lines $x$ and $y$.
We show that the point $O=x\cap y$ lies on the line $\{O_{T_t, T_{\lambda}} \}$ for all
$\lambda \neq 0, 1$, except for the case $\lambda = h$,
in which $H$ is the orthocenter of $A_hB_hC_h$. Moreover, we show that

\otstup

{\it $O$ is the common orthology center for a given triangle $T_{\lambda}$
($\lambda \neq 0, 1, h$) and another triangle $T_{\mu}$ from our family,
with the correspondence $\lambda \leftrightarrow \mu$ being linear fractional.}

\otstup

{\it Proof.} Let $a\cap x = (x_a, 0)$, $a\cap y = (0, y_a)$,
and so on. Then the coordinates of the vertices of the triangle $T_{\lambda}$ are:
$A_{\lambda}=(\lambda x_a, (1-\lambda ) y_a)$, and so on.
The condition $OA_{\mu} \perp B_{\lambda}C_{\lambda}$ can be written as
$$\mu \lambda x_a (x_c-x_b) + (1-\mu)(1-\lambda) y_a (y_c-y_b)=0.$$
As we can see, this equation is symmetric with respect to $\lambda$ and $\mu$, i.e.,
it is equivalent to the perpendicularity $OA_{\lambda} \perp B_{\mu}C_{\mu}$.
If 
$$\lambda \neq \frac{1}{k_a+1}, 
\ \ \text{where}  \ \ 
k_a = \frac{x_a (x_c-x_b)}{y_a (y_c-y_b)},$$ then this equality transforms to $$\mu = \frac{1-\lambda}{1-\lambda(k_a+1)}.$$ The calculation shows that
the condition of equality of the corresponding coefficients $k_a=k_b$ (or $k_a=k_c$) is equivalent to
the concurrence of the lines $a$, $b$, $c$ (whose equations are $x/x_a+y/y_a=1$, etc.).
\edok

\otstup 

It is clear that the singular value $\lambda = 1 / (k_a+1)$ corresponds to the fact that $H$ is the orthocenter of the triangle $T_{\lambda}$.

In the correspondence $\lambda \leftrightarrow \mu$ obtained above, the equality $\lambda = \mu$ is possible for no more than two values of $\lambda$.
This equality corresponds to the fact that $O$ serves as the orthocenter for $T_{\lambda}$.

\otstup

Note that the coincidence of orthology centers $O_{T_{\mu}, T_{\lambda}} = O_{T_{\lambda}, T_{\mu}}$ is in principle possible in no more than two points of intersection of $\{O_{T_{t}, T_{\lambda}}\}$ and the hyperbola $\{O_{T_{\lambda}, T_{t}} \}$. One of these points is $O$, and as observed, at $\lambda \neq 0, 1, h$, we indeed have the coincidence of orthology centers, and the other point is the orthocenter of the triangle $T_{\lambda}$.

\otstup

{\bf The case of degenerate $T_{\infty}$}


Consider the case of the degenerate triangle $T_{\infty}$. Let $T_0$ be the second degenerate triangle lying on the line $y$. We show that in this case, for $\lambda \neq 0$,
\begin{center} 
{\it
the line $\{O_{T_{t} , T_{\lambda}}\}$ is parallel to the line $y$.}
\end{center}

{\it Proof.} Straightforward calculation in the affine coordinate system $Oxy$ from \ref{1555} shows that the abscissa of the point $O_{T_{\mu} , T_{\lambda}}$ (i.e., the point of intersection of the perpendiculars from $A_{\mu}$ to $B_{\lambda}C_{\lambda}$ and from $B_{\mu}$ to $A_{\lambda}C_{\lambda}$) does not depend on $\mu$.
\edok

\otstup

We see that in this case, the hyperbola $\{O_{T_{\lambda}, T_{t}}\}$ and the line $\{O_{T_{t} , T_{\lambda}}\}$ have only one common point (the orthocenter of the triangle $T_{\lambda}$).
\otstup

{\bf Example: the problem of flies on gergonnians}

{\it Problem.} (E. Bakaev)
The inscribed circle of triangle $ABC$ touches its sides at points $A_1, B_1, C_1$. Three flies crawled along the lines $AA_1, BB_1, CC_1$ with constant speeds so that at one moment they were at points $A, B, C$, and at another moment they were at points $A_1, B_1, C_1$.
Prove that at the moments when the flies were on one straight line, the center of the inscribed circle was also on this line, and there are two such moments, with the corresponding lines being perpendicular.
(The lines $AA_1, BB_1, CC_1$ are called gergonnians of the triangle, they intersect at the Gergonne point.) 

\otstup 
A weakened version of this problem (in which it was already given that there were two moments of collinearity) was proposed at the Kolmogorov Cup in 2017.

{\it Solution.} In the problem, we are dealing with a perspective orthologic family generated by triangles $ABC$ and $A_1B_1C_1$. Note that $ABC$ and $A_1B_1C_1$ are orthologic, and both centers of orthology coincide with $I$.
From this (or by directly showing that $(\vek{BB_1} - \vek{AA_1}) \nparallel
(\vek{CC_1} - \vek{AA_1})$), we see that
this is not the case of degeneration $T_{\infty}$.
Based on what was proven above in \ref{371000}, we understand that $I$ is the point of intersection of the perpendicular lines $x$ and $y$, 
on which the degenerate triangles of this family lie.
\edok

\bigskip 
\bigskip 
\bigskip 

\section{Interpretation in 
\texorpdfstring{$\mathbb{R}^4$}{ h} 
and 
\texorpdfstring{$\mathbb{R}P^3$}{ h}
}

We represented the triangle $ABC$, up to translation, by an ordered pair of Euclidean vectors $(\vek{b}, \vek{c}) = (\vek{AB}, \vek{AC})$. This pair $(\vek{b}, \vek{c})$ can be identified with a vector in $\mathbb{R}^4 = \mathbb{R}^2\times \mathbb{R}^2$.

Next, we discuss how some of the previous scenarios transfer to this model and what they mean in the language of linear algebra.

\begin{predl}\label{51}
{\bf Classes of homothets}
\end{predl}

The equivalence class $[ABC]$ can also be uniquely determined by pairs of the form $(\alpha \vek{AB}, \alpha \vek{AC})$, such that $[ABC]$ corresponds to a one-dimensional subspace in the vector space $\mathbb{R}^4$. Therefore, the set of classes (quotient set) can now be identified with the projective space $\mathbb{R}P^3$.


\begin{predl}\label{52}
{\bf Linear families}
\end{predl}

From formulas (\ref{linlin}), (\ref{linob}), we see that a linear family of triangles $\{A_tB_tC_t\}$ corresponds to a line (i.e., a one-dimensional linear submanifold) in $\mathbb{R}^4 = \mathbb{R}^2\times \mathbb{R}^2$, and a linear family of classes (containing at least two distinct classes) corresponds to a two-dimensional subspace in $\mathbb{R}^4 = \mathbb{R}^2\times \mathbb{R}^2$, which is the linear span of vectors $(\vek{b}_0, \vek{c}_0)$ and $(\vek{b}_1, \vek{c}_1)$.


\begin{predl}\label{53}
{\bf Linear family of classes as a graph}
\end{predl}

Equation (\ref{linob}) can be written in the following form:
\begin{equation}\label{phi}
(\vek{b}, \vek{c} ) = (\vek{v}, \varphi (\vek{v})),
\end{equation}
where $\varphi \,:\, \mathbb{R}^2 \to \mathbb{R}^2 $ is a linear bijective operator (automorphism) on $\mathbb{R}^2$, defined by its action on basis vectors: 
$\varphi (\vek{b}_t) = \vek{c}_t$, $t=0, 1$ (here, as before, we assume $\vek{b}_t = \vek{A_{t}B_{t}}$, $\vek{c}_t = \vek{A_{t}C_{t}}$). For instance, in Example 1 from \ref{17}, the corresponding operator $\varphi$ is a spiral similarity.

The fact that $(\vek{b}_0, \vek{b}_1)$ and $(\vek{c}_0, \vek{c}_1)$ are indeed bases is equivalent to the pairs of moving points $A, B$ and $A, C$ being nonsingular.

The condition that the pair $B, C$ is singular means that $\vek{B_0C_0}\parallel \vek{B_1C_1}$, which is equivalent  $\varphi (\vek{b}_0) - \vek{b}_0 \parallel \varphi (\vek{b}_1) - \vek{b}_1$, which in turn is equivalent to the degeneracy of the operator $\varphi - \Id$,
or the fact that $\lambda = 1$ is an eigenvalue of $\varphi$.
Conversely, according to (\ref{phi}), a bijective linear operator $\varphi$, for which $\lambda = 1$ is not an eigenvalue, uniquely determines a nonsingular linear family of classes of triangles.


\begin{predl}\label{54}
{\bf Degenerate triangles}
\end{predl}

Let us consider the coordinates $(x_1, x_2, x_3, x_4)$, where $\vek{b} = (x_1, x_2)$, $\vek{c} = (x_3, x_4)$. The degeneracy condition for the corresponding triangle $ABC$ is given by $x_1x_4-x_2x_3=0$. We observe that the degenerate triangles $ABC$ correspond under the mapping $ABC \,\mapsto \, (\vek{b}, \vek{c})$ to the asymptotic cone of a quadratic form with the signature $(2,2)$.

The  condition of degeneracy of a triangle can also be reformulated differently, based on (\ref{phi}). We see that a degenerate triangle corresponds to an eigenvector of the operator $\varphi$ introduced earlier (and the class of homothets of a degenerate triangle corresponds to a one-dimensional invariant subspace of the operator~$\varphi$).

As we noted in \ref{15}, in a linear family of triangles, there are either no more than two degenerate triangles, or all triangles of the family are degenerate. Accordingly, in a linear family of homothets, there are either no more than two classes of degenerate triangles, or all classes consist of degenerate triangles. Now we can prove it in a different way, based on $\varphi$: either $\varphi$ is a homothety (and then the corresponding class consists only of degenerate triangles), or $\varphi$ has no more than two eigenvectors, up to proportionality.


\begin{predl}\label{55}
{\bf Orthology and the self-adjoint operator}
\end{predl}

Next, we derive a condition on the operator $\varphi$ from (\ref{phi}) that is equivalent to the orthology of the linear family of classes.

The orthology condition for any two triangles $ABC$, $A'B'C'$ from our family: $\vek{b} \vek{c'} - \vek{b'} \vek{c}=0$ implies the equality $ \vek{b} \varphi (\vek{b'}) - \vek{b'} \varphi (\vek{b})=0$
for all $\vek{b}$ and $\vek{b'}$, i.e., it is equivalent to $\varphi $ being a self-adjoint operator. Thus, precisely the graphs of self-adjoint operators $\varphi : \mathbb{R}^2 \to  \mathbb{R}^2 $ correspond to orthologic families of classes.

\otstup

The obtained correspondence now enables us to prove the statement from \ref{34} 
about the degenerate triangles in an orthologic family as follows.
A well-known linear algebra theorem states:
a linear operator in an $n$-dimensional Euclidean space is self-adjoint 
if and only if it has an orthogonal basis composed of its eigenvectors. 
In our case where $n=2$, the operator $\varphi:\, \mathbb{R}^2 \to  \mathbb{R}^2$ is self-adjoint 
if and only if it has a pair of orthogonal eigenvectors.

\begin{predl}\label{56}
{\bf Orthology as skew-orthogonality in $\mathbb{R}^4$}
\end{predl}

As we recall (see \ref{25}), the orthology of a triangle $ABC$ can be interpreted as 
a relation in $\mathbb{R}^4=\mathbb{R}^2\times \mathbb{R}^2$.

We define a multiplication $*$ on $\mathbb{R}^2\times \mathbb{R}^2 = \{(\vek{b}, \vek{c})\}$ 
by the equality 
\begin{equation}\label{kosoR^4}
(\vek{b}, \vek{c}) * (\vek{b}', \vek{c}') = \vek{b} \vek{c}' - \vek{b}' \vek{c}.
\end{equation}
This multiplication $*$ establishes a bilinear skew-symmetric nondegenerate form on $\mathbb{R}^4$. 
In accordance with (\ref{ort_pary}), the orthology relation 
corresponds to <<skew-orthogonality>> with respect to $*$:
\begin{equation}\label{definR^4}
ABC\ort A'B'C' \,\,\, \Longleftrightarrow \,\,\, 
(\vek{b}, \vek{c}) * (\vek{b}', \vek{c}') = 0.
\end{equation}


\begin{predl}\label{57}
{\bf Skew-orthogonal complements}
\end{predl}

We denote the orthogonal complement to a subspace $U$ 
with respect to the introduced bilinear multiplication $*$ (or skew-orthogonal complement) by $U^{\bot}$.

Now, in $\mathbb{R}^4$, the set of triangles that are orthologic to a given triangle $ABC$ can be described as follows:
it is the (three-dimensional) skew-orthogonal complement to the (one-dimensional) subspace generated by $(\vek{b}, \vek{c})$, that is, the subspace $\lin{ (\vek{b}, \vek{c})}^{\bot}$.

The set of triangles that are orthologic to both $ABC$ and $A'B'C'$ from different classes 
$[ABC]$ and $[A'B'C']$ ---
is the (two-dimensional) subspace $\lin{ (\vek{b}, \vek{c}), (\vek{b}', \vek{c}')}^{\bot}$.

The case $ABC\ort A'B'C'$ 
is characterized by the restriction of the multiplication $*$ to the subspace $\lin{ (\vek{b}, \vek{c}), (\vek{b}', \vek{c}')}$ being zero,
which means the subspace 
$\lin{ (\vek{b}, \vek{c}), (\vek{b}', \vek{c}')}$ coincides with its skew-orthogonal complement.
In other words, in this case, a triangle orthologic to both $ABC$ and $A'B'C'$ 
necessarily belongs to the linear family of classes generated by the classes $[ABC]$ and $[A'B'C']$.

Subspaces with the property $U^{\bot}=U$ are called {\it Lagrangian} (with respect to $*$).
Thus, in the language of orthology, we have come to the known correspondence
between the graphs of self-adjoint operators and Lagrangian subspaces.

\end{document}